\documentclass{gtpart}
\usepackage{pinlabel}
\usepackage{amssymb,amsfonts}
\usepackage[utf8]{inputenc}
\usepackage{amsmath}
\usepackage[all,arc]{xy}
\usepackage{enumerate}
\usepackage{mathrsfs}
\usepackage{tikz}
\usepackage[utf8]{inputenc}
\usetikzlibrary{matrix,arrows}
\allowdisplaybreaks
\title[mmmmmm]{jjjjjjjjjjjjjjjjjjjjjjj}
\title{Samplehjhhhhhr layout}
\author[Hassan Jolany]{Hassan Jolany}
\givenname{Hassan }
\surname{Jolany}
\address{Laboratoire de Mathématiques Paul Painlevé,\\\newline 
Université des Sciences et Technologies de Lille,\\\newline
Lille\\Fance}
\email{hassan.jolany@math.univ-lille1.fr}
\subject{primary}{}{53C44}
\arxivreference{math.AW/0611679}
\arxivpassword{5spud}
\volumenumber{}
\issuenumber{}
\publicationyear{}
\papernumber{}
\startpage{}
\endpage{}
\doi{}
\MR{}
\Zbl{}
\received{}
\revised{}
\accepted{}
\published{}
\publishedonline{}
\proposed{}
\seconded{}
\corresponding{}
\editor{}
\version{}
\newtheorem{thm}{Theorem}[section]

\theoremstyle{definition}
\newtheorem{defn}[thm]{Definition}

\numberwithin{equation}{section}
\makeop{Homo}

\theoremstyle{definition}

\usepackage{pinlabel}  
\usepackage{graphicx}  
\usepackage[all]{xy}
\title{Canonical metric on moduli spaces of log Calabi-Yau varieties}

\author{}
\givenname{nnnn}
\surname{uuuuuu}
\address{}
\email{}
\urladdr{}


\arxivreference{}
\arxivpassword{}

\volumenumber{}
\issuenumber{}
\publicationyear{}
\papernumber{}
\startpage{}
\endpage{}
\doi{}
\MR{}
\Zbl{}
\received{}
\revised{}
\accepted{}
\published{}
\publishedonline{}
\proposed{}
\seconded{}
\corresponding{}
\editor{}
\version{}

\begin{document}
\begin{abstract} 

In this paper, by applying Greene-Shapere-Vafa-Yau semi-flat metric, we give a new proof of closed formula of Weil-Petersson metric on moduli space of Calabi-Yau varieties.

\end{abstract}
\maketitle{}
\section*{Introduction}

In this note we try to find the K\"ahler potential of logarithmic Weil-Petersson metric on moduli space of log Calabi-Yau varieties. We use the analysis of semi Ricci-flat metric introduced in \cite{6}. 
Historically, A. Weil introduced a K\"ahler metric for the Teichmuller space $\mathcal T_{g,n}$, the space of homotopy marked Riemann surfaces of genus $g$ with $n$ punctures and negative Euler characteristic. The Weil-Petersson metric measures the variations of the complex structure of ${R}$. The quotient of the Teichmuller space $\mathcal T_{g,n}$ by the action of the mapping class group is the moduli space of Riemann surfaces $\mathcal M_{g,n}$. The Weil-Petersson metric is mapping class group invariant and descends to $\mathcal M_{g,n}$.  A.Fujiki \cite{24} considered Weil-Petersson metric on moduli space of general type K\"ahler-Einstein varieties. In this note we consider the logarithmic Weil-Petersson metric on moduli space of log Calabi-Yau varieties(if exists!).

We start we some elementary definitions of relative K\"ahler metric from Mourougane and Takayama \cite{29}.

\begin{defn}

Let $\pi:X\to Y$ be a holomorphic map of complex manifolds. A real d-closed $(1,1)$-form $\omega$
on $X$ is said to be a relative K\"ahler form
for $\pi$, if for every point $y\in Y$
, there exists an
open neighbourhood
$W$
of
$y$
and a smooth plurisubharmonic function
$\Psi$
on
$W$ such that $\omega+\pi^*(\sqrt[]{-1}\partial\bar\partial\Psi)$ is a K\"ahler form on $\pi^{-1}(W)$. A morphism
$\pi$ is said to be K\"ahler, if there exists a relative K\"ahler form for $\pi$, and $\pi:X\to Y$
is said to be a K\"ahler fiber space, if $\pi$ is proper, K\"ahler, and surjective with connected fibers.
\end{defn}

We consider an effective holomorphic family of complex manifolds. This means we have a holomorphic map $\pi:X\to Y$
between complex manifolds such that
\;

1.The rank of the Jacobian of $\pi$ is equal to the dimension of $Y$ everywhere.

2.The fiber $X_t=\pi^{-1}(t)$ is connected for each $t\in Y$

3.$X_t$ is not biholomorphic to $X_{t'}$ for distinct points
$t;t'\in B$.

It is worth to mention that Kodaira showed that all fibers are dieomorphic to each other.

For any effective holomorphic family of compact manifolds
$\pi:X \to Y$ of dimension $n$ with fibers
$X_y$
for
$y\in Y$
the Calabi-Yau forms $\omega_{X/Y}$ depend differentiably on the parameter $y$. The relative K\"ahler form is denoted by

$$\omega_{X/Y}=\sqrt[]{-1}g_{\alpha,\bar\beta}(z,y)dz^\alpha\wedge d\bar z^\beta$$ Moreover take $\omega_X=\sqrt[]{-1}\partial\bar\partial\log \det g_{\alpha,\bar\beta}(z,y)$ on the total space $X$. The fact is $\omega_X$ in general is not K\"ahler on total space and $\omega_X|_{X_y}=\omega_{X_y}$. More precisely $\omega_X=\omega_F+\omega_H$ where $\omega_F$ is a form along fiber direction and $\omega_H$ is a form along horizontal direction. $\omega_H$
may not be K\"ahler metric in general, but $\omega_F$ is K\"ahler metric.
Now let $\omega$ be a relative K\"ahler form on $X$ and $m:=\dim X-\dim Y$, We define the relative Ricci form $Ric_{X/Y,\omega}$ of $\omega$ by

$$Ric_{X/Y,\omega}=-\sqrt[]{-1}\partial\bar\partial\log (\omega^m\wedge \pi^*|dy_1\wedge dy_2\wedge...\wedge dy_k|^2)$$
where $(y_1,...,y_k)$ is a local coordinate of $Y$. Here $Y$ assumed to be a curve

Let for family $\pi:\mathcal X\to Y$
$$\rho_{y_0}:T_{y_0}Y\to H^1(X,TX)=\mathcal H_{\bar\sigma}^{0,1}(TX)$$ be the Kodaira–Spencer map for the corresponding deformation of $X$ over $Y$ at the
point $y_0\in Y$ where $\mathcal X_{y_0}=X$

If $v\in T_{y_0}Y$ is a tangent vector, say $v=\frac{\partial}{\partial y}\mid_{y_0}$ and $\frac{\partial}{\partial s}+b^\alpha\frac{\partial}{\partial z^\alpha}$ is any lift
to $\mathcal X$ along $X$, then 

$$\bar\partial \left(\frac{\partial}{\partial s}+b^\alpha\frac{\partial}{\partial z^\alpha}\right)=\frac{\partial b^\alpha(z)}{\partial z^{\bar\beta}}\frac{\partial}{\partial z^\alpha}dz^{\bar\beta}$$
is a $\bar\partial$-closed form on $X$, which represents $\rho_{y_0}(\partial/\partial y)$. 

The Kodaira-Spencer map is induced as edge homomorphism by the short exact sequence
$$0
\to T_{X/Y}\to TX
\to  \pi^*T_Y\to  0$$

We briefly explain about the Weil-Petersson metric on  moduli space of polarized Calabi-Yau manifolds. We study the moduli space of Calabi-Yau manifolds via the Weil-Petersson metric. We outline the imortant properties of such metrics here.

The Weil-Petersson metric is not complete metric in general but in the case of abelian varieties and $K3$ surfaces, the Weil-Petersson
metric turns out to be equal to the Bergman metric of the Hermitian symmetric
period domain, hence is in fact complete K\"ahler Einstein metric. Weil  and Ahlfors  showed that the Weil-Petersson metric is a K\"ahler metric. Ahlfors proved that it has negative holomorphic sectional, scalar, and Ricci curvatures. The quasi-projectivity of coarse moduli spaces of polarized Calabi-Yau manifolds in the category of separated analytic spaces (which also can be constructed  in
the category of Moishezon spaces) has been proved by Viehweg\cite{33}. By using Bogomolov-Tian-Todorov theorem, these moduli spaces are smooth K\"ahler orbifolds equipped with the Weil-Petersson metrics. Let $X \to M$ be a family of polarized Calabi-Yau manifolds. Lu and Sun showed that the volume of the first Chern class with respect to the Weil-Petersson metric over the moduli space M is a rational number. Gang Tian (also Georg Schumacher) proved that the Weil-Petersson metric  on  moduli space of polarized Calabi-Yau manifolds is just pull back of Chern form of the tautological of $\mathbb CP^N$ restricted to period domain which is an open set of a quadric in $\mathbb CP^N$ and he showed that holomorphic sectional curvature is bounded away from zero. Let $X$ be a compact projective Calabi-Yau manifold and let $f:X\to Y$ be an algebraic fiber space with $Y$ an irreducible normal algebraic variety of lower dimension then Weil-Petersson metric measuring the change of complex structures of the fibers.

Now, consider a polarized K\"ahler manifolds $\mathcal X\to S$ with K\"ahler metrics $g(s)$ on $\mathcal X_s$. We can define a possibly degenerate hermitian metric $G$ on $S$ as follows:
Take Kodaira-Spencer map $$\rho:T_{S,s}\to H^1(X,T_X)\cong H^{0,1}_{\bar{\partial}}(T_X)$$
into harmonic forms with respect to
$g(s)$; so for $v,w\in T_s(S)$ , we
may define 

$$G(v,w):=\int_{\mathcal X_s}<\rho(v),\rho(w)>_{g(s)}$$

When $\mathcal X \to S$ is a polarized K\"ahler-Einstein family and $\rho$ is injective $G_{WP}:=G$ is called the Weil-Petersson metric on $S$. Tian-Todorov, showed that if we take $\pi:\chi\to S$, $\pi^{-1}(0)=X_0=X$, $\pi^{-1}(s)=X_s$ be the family of $X$, then $S$ is a non-singular complex analytic space such that $$dim_{\mathbb C}S=dim_{\mathbb C}H^1(X_s,TX_s)$$

Note that in general, 
if
$f:
X
\to S$
be a smooth projective family over a complex manifold
$S$. Then for every positive integer
$m$, $$P_m(X_s)=dim H^0(X_s,\mathcal O_{X_s}(mK_{X_s}))$$
is locally constant function on $S$. 

It is worth to mention that the fibers $X_s$ are diffeomorphic to each other and if fibers $X_s$ be biholomorphic then $\pi$ is holomorphic fiber bundle and Weil-Petersson metric is zero in this case in other words the Kodaira-Spencer maps
$$\rho:T_{S,s}\to H^1(X_s,T_{X_s})\cong H^{0,1}_{\bar{\partial}}(T_{X_s})$$
are zero. In special case, let $dim X_s=1$, then the fibers are elliptic curves and $\pi$ is holomorphic fiber bundle and hence the Weil-Petersson metric is zero. In general, the Weil-Petersson metric is semipositive definite on the moduli space of Calabi-Yau varieties. Note that Moduli space of varieties of general type has Weil-Petersson metric. The moduli space of K-stable varieties admit Weil-Petersson metric also.\cite{37},\cite{38}

\;
One of important thing which we need to say, is that we assume the central fiber of our degeneration is smoothable and in fact note that not any variety can be smoothable, for example there exists a Calabi-Yau variety(due to Namikawa) which remain singular under any flat deformation. Even when the central fiber has nice variety like Calabi-Yau variety with mild singularity in the sense of MMP, we have a few result for existence of smooth degeneration. Let $\pi : X \to \Delta$ be a proper map of a variety $X$ onto $\Delta$ such that

1. $X-\pi^{-1}(0)$ is smooth.

2. The fibers $X_t = \pi^{-1}(t)$ are nonsingular compact K\"ahler manifolds for every $t\not= 0$

We call $\pi$ a degeneration of $X_t$ ( $t \not= 0$), the fiber $X_0 = \pi^{-1}(0)$ the degenerate fiber
and $X$ the total space. We say that $X_t$ ( $t \not= 0$ ) has a degeneration to a variety
$X_0$ and that $X_0$ is smoothable to $X_t$ or has a smoothing to $X_t$
. We say that $X$ is a
trivial degeneration if $X_0$ is smooth.

Friedman\cite{73} showed that any Semistable(we explain this notion later) $K3$ surface can be smoothed into a smooth $K3$ surface under a flat deformation with a smooth total space.Y. Namikawa\cite{74} showed the following result. Let $X$ be a Calabi–Yau threefold with terminal singularities.

(1) If $X$ is $\mathbb Q$-factorial; then $X$ is smoothable.

(2) If every singularity of $X$ is different from an ordinary double point; then $X$ is smoothable.

A reduced complex analytic space $X$ of dimension $n$ is a
normal crossing variety (or n.c.variety) if for each point $p \in X$,
$$\mathcal O_{X,p}\cong \frac{\mathbb C\{x_0,x_1,...,x_n\}}{(x_0x_1...x_r)}\;\;\;\; 0\leq r=r(p)\leq n $$.

In addition, if every component $X_i$ of $X$ is smooth, then $X$ is called a, simple normal crossing variety (or s.n.c.variety).

Let $D=\text{Sing}(X)$, and $X_i$ be a component of $X$ and let $I_X$, (resp. $I_D$) be the defining ideal
of $X_i$ (resp. $D$) in $X$. Then define

$$\mathcal O_D(-X)=I_{X_1}/I_{X_1}I_D\otimes_{\mathcal O_D}...\otimes_{\mathcal O_D} I_{X_m}/I_{X_m}I_D$$ and take $\mathcal O_D(X):=\mathcal O_D(-X)^\vee$

A normal crossing variety $X$ is called $d$-semistable if its infinitesmal normal bundle $\mathcal O_D(X)$ be trivial.

Kawamata and Namikawa \cite{75} proved the following theorem for normal crossing varieties with some additional assumption to get smoothing of a flat degeneration

Let $X$ be compact K\"ahler d-semi-stabte n.c.variety of dimension $n\geq 3$ and  let $X^{[0]}$ be the normalization of $X$. Assume the following conditions:

 (a) $\omega_X\cong \mathcal O_X$,

 (b) $H^{n-1}(X, \mathcal O_X) = 0$, and

 (c) $H^{n-2}(X^{[0]},\mathcal O_{X^{[0]}}) = 0$

Then $X$ is smoothable by a flat deformation.

\;

Let us explain briefly how is organized this paper. In Section, Introduction, we provide a brief introduction on deformation theory of K\"ahler metrics, geometry of Weil-Petersson metrics and . In prove the main theorem of this paper about explicit formula for canonical metric on moduli space of log Calabi-Yau varieties by using foberwise Calabi-Yau metric and explain Bismut-Gillet-Soule fiberwise integral formula and its connection with Weil-Petersson metric and fiberwise K\"ahler-Einstein metric. We give a relation between relative volume form constructed by fiberwise Calabi-Yau metric and logarithmic Weil-Petersson metric. Moreover we talk about fiberwise Calabi-Yau foliation and extend such foliation for logarithmic case. In last part of this paper we talk about A relation between Invariance of plurigenera and positivity of logarithmic Weil-Petersson metric by using relative K\"ahler Ricci flow method (introduced by R.Berman) and also we show that Song-Tian-Tsuji measure is bounded if and only if the central fiber has log terminal singularities at worst. 

\;
\;

\textbf{Remark}: Let $(E,\|.\|)$ be the direct image bundle $f_*(K_{X'/S})$, where $X'=X\setminus D$, of relative canonical line bundle equipped with the $L^2$ metric $\|.\|$. Then the fibre $E_y$ is $H^0(X_y\setminus D_y,K_{X_y\setminus D_y})$. Since the pair $(X_y,D_y)$ is Calabi-Yau pair, hence $H^0(X_y\setminus D_y,K_{X_y\setminus D_y})$ is a 1-dimensional vector space. This implies that $E$ is a line bundle.

We give a new proof to the following theorem \cite{18},\cite{28}.
\begin{thm}
Take holomorphic fiber space $\pi:X\to B$ and assume $\Psi_y$ be any local non-vanishing holomorphic section of Hermitian line bundle $\pi_*(K_{X/B}^l)$, then the Weil-Petersson (1,1)-form on a small ball $N_r(y)\subset B$ can be written as

$$\omega_{WP}=-\sqrt{-1}\partial_y\bar{\partial_y}\log \left((\sqrt{-1})^{n^2}\int_{X_y}(\Psi_y\wedge \overline{\Psi_y})^{\frac{1}{l}}\right)$$
Note that $\omega_{WP}$ is globally defined on $B$

\end{thm}

Now because we are in deal with Calabi-Yau pair $(X,D)$ which $K_X+D$ is numerically trivial so we must introduce Log Weil-Petersson metrics instead Weil-Petersson metric. Here we introduce such metrics on moduli space of paired Calabi-Yau fibers $(X_y,D_y)$. Let $i:D\hookrightarrow X$ and $f : X \to Y$ be holomorphic mappings of complex manifolds
such that $i$ is a closed embedding and $f$ as well as$ f \circ i$ are proper and smooth. Then a holomorphic family $(X_y,D_y)$ are the fibers $X_y=f^{-1}(y)$ and $ D_y=(f \circ i)^{-1}(y)$. Such family give rise to a fibered groupoid $p:\bf F\to \bf A$ from of category $\bf F$ to the category of complex spaces with distinguished point in the sense of Grothendieck, \cite{30}. There exists the moduli space of $\mathcal M$ of such family because any $(X_y,D_y)$ with trivial canonical bundle is non-uniruled. Now $X\setminus D$ is quasi-projective so we must deal with quasi-coordinate system instead of coordinate system. Let $(X,D)$ be a Calabi-Yau pair and take $X'=X\setminus D$ equipped with quasi-coordinate system. We say that a tensor $A$ on $X'$ which are covariant of type $(p,q)$ is quasi-$C^{k,\lambda}$-tensor, if it is of class $C^{k,\lambda}$ with respect to quasi-coordinates. Now we construct the logarithmic version of Weil-Petersson metric on moduli space of paired Calabi-Yau fibers $f:(X,D)\to Y$.

Now, because we are in deal with singularities, so we use of $(1,1)$-current instead of $(1,1)$-forms which is singular version of forms. A current is a differential form with distribution coefficients. Let, give a definition of current here.
We recall a singular metric $h_{\text{sing}}$ on a Line bundle $L$ which locally can be written as $h_{\text{sing}}=e^\phi h$ where $h$ is a smooth metric, and $\phi$ is an integrable function. Then one can define locally the closed current $T_{L,h_{\text{sing}}}$ by the following formula

$$T_{L,h_{\text{sing}}}=\omega_{L,h}+\frac{1}{2i\pi}\partial\bar{\partial}\log \phi$$

The current Geometry is more complicated than symplectic geometry. For instance in general one can not perform the wedge product of currents due to this fact that one can not multiply the coefficients which are distributions and in general the product of two distributions is not well defined. However, in some special cases one can define the product of two currents. Here we mention the following important theorem about wedge product of two currents

One can simply defines the space of currents
to be the dual of space of smooth forms, defined as forms on the regular part
$X_{
reg}$
which, near
$X_
{sing}$
, locally extend as smooth forms on an open set of
$\mathbb C^
N$
in which
$X$
is locally embedded. A K\"ahler current on a compact complex space
$X$
is a
closed positive current
$T$
of bidegree
$(1,1)$ which satisfies
$T
\geq
\epsilon\omega$
for some
$\epsilon>
0 $ and some smooth positive hermitian form
$\omega$
on
$X$. In fact, This is a real closed current of type $(1,1)$, that is a linear form on the space of compactly supported forms of degree $2n-2$ on $X$, and $n=dim X$. Mre precisely, Let $\mathcal A^{p,q}_c(X)$ denote the space of $C^{\infty}{(p,q)}$ forms of compact support on $X$ with usual Fréchet space structure. The dual space $\mathcal D^{p,q}(X):=\mathcal A_c^{n-p,n-q}(X)^*$ is called the space of $(p,q)-$ currents on $M$. The Linear operators $\partial:\mathcal D^{p,q}(X)\to \mathcal D^{p+1,q}(X)$ and $\partial:\mathcal D^{p,q}(X)\to \mathcal D^{p,q+1}(X)$ is defined by 

$$\partial T(\varphi)=(-1)^{p+q+1}T(\partial \varphi), \; T\in \mathcal D^{p,q}(X), \; \varphi\in A_c^{n-p-1,n-q}(X)$$

and

$$\bar\partial T(\varphi)=(-1)^{p+q+1}T(\bar\partial \varphi), \; T\in \mathcal D^{p,q}(X), \; \varphi\in A_c^{n-p,n-q-1}(X)$$
We set $d=\partial+\bar\partial$. $T\in \mathcal D^{p,q}(X)$ is called closed if $d T=0$. $T\in \mathcal D^{p,p}(X)$ is called real if $T(\varphi)=T(\bar \varphi)$ holds for all $\mathcal A^{n-p,n-q}_c(X)$. A real $(p,p)$-current $T$ is called positive if $(\sqrt{-1})^{p(n-p)}T(\eta\wedge\bar\eta)\geq 0$ holds for all $\eta\in \mathcal A_c^{p,0}(X)$. 

The topology on space of currents are so important. In fact the space of currents with weak topology is a
Montel space, i.e., barrelled,
locally convex, all bounded subsets are precompact which here barrelled topological vector space is Hausdorff topological vector space for which every barrelled set in the space is a neighbourhood for the zero vector.

\;

Also because we use of push-forward and Pull back of a current and they can cont be defined in sense of forms, we need to introduce them. If $f:X\to Y$ be a holomorphic map between two compact K\"ahler manifolds then one can push-forward a current $\omega$ on $X$ by duality setting 

$$\left<f_*\omega,\eta \right>:=\left<\omega,f^*\eta\right>$$

In general, given a current $T$ on $Y$, it is not possible to define its pull-back by a holomorphic map. But it is possible to define pull-back of positive closed currents of bidegree $(1,1)$. We can writes such currents as $T=\theta+dd^c\varphi$ where $\theta\in {T}$ is a smooth form, and thus one define the pull-back of current $T$ as follows

$$f^*T:=f^*\theta+dd^c\varphi \circ f$$

Let $X$ and $Y$ be compact K\"ahler manifolds and let $f:X\to Y$ be the blow up of $Y$ with smooth connected center $Z$ and $\omega\in H^{1,1}(X,\mathbb R)$. Demailly showed that

$$\omega=f^*f_*\omega+\lambda{E}$$ where
$E$ is the exceptional divisor and $\lambda\geq -v(\omega,Z)$ where $v(\omega,Z)=\inf_{x\in Z}v(\omega,x)$ and $v(\omega,x)$ is the Lelong number.

Note that if $\omega$ is K\"ahler then,

$$d\text{Vol}_{\omega_y}(X_y)=df_*(\omega^n)=f_*(d\omega^n)=0$$
So, $\text{Vol}(X_y)=C$ for some constant $C>0$ for every $y\in Y$ where $\pi^{-1}(y)=X_y$ is general smooth fiber. Moreover direct image of volume form $f_*\omega_X^n=\sigma\omega_{Y}^m$ where $\sigma\in L^{1+\epsilon}$ for some positive constant $\epsilon$

\begin{thm}

If $T$ is a positive $(1,1)$-current then locally one can find a plurisubharmonic
function
$u$
such that
$$\sqrt{-1}\partial\bar\partial u=T$$
Note that, if $X$ be compact then there is no global plurisubharmonic function $u$\cite{31}.

\end{thm}

\textbf{Lelong number}:
Let
$W\subset \mathbb C^n$
be a domain, and $\Theta$ a positive current of degree $(q,q)$ on
$W$. For a point $p\in W$
one defines
$$\mathfrak v(\Theta,p,r)=\frac{1}{r^{2(n-q)}}\int_{|z-p|<r}\Theta(z)\wedge (dd^c|z|^2)^{n-q}$$
The
Lelong number
of $\Theta$ at
$p$
is defined as

$$\mathfrak v(\Theta,p)=\lim_{r \to 0}\mathfrak v(\Theta,p,r)$$

Let $\Theta$ be the curvature of singular hermitian metric $h=e^{-u}$, one has

$$\mathfrak v(\Theta,p)=\sup\{\lambda\geq 0: u\leq \lambda\log(|z-p|^2)+O(1)\}$$
see \cite{36}

\section*{Main Theorem}

Now we are ready to state our theorem. We must mention that The result of Tian, and Candelas et al. was on Polarized Calabi-Yau fibers and in this theorem we consider non-polarized fibers.

\begin{thm}
Let $\pi:X\to Y$ be a smooth family of compact K\"ahler manifolds whith Calabi-Yau fibers. Then Weil-Petersson metric can be written as $$\omega_{WP}=-\sqrt[]{-1}\partial_y\bar\partial_y\log\int_{X_y}|\Omega_y|^2$$
where  $\Omega_y$ is a holomorphic $(n,0)$-form
on $\pi^{-1}(U)$, where
$U$ is a neighborhood of y
\end{thm}
{Proof}: For prrof, We need to recall the Yau-Vafa semi Ricci flat metrics. Since fibers are Calabi-Yau varieties, so $c_1(X_y)=0$, hence there is a smooth function $F_y$ such that $Ric(\omega_y)=\sqrt{-1}\partial\bar{\partial}F_y$ . The function $F_y$ vary smoothly in $y$. By Yau's theorem, there is a unique Ricci-flat K\"ahler metric $\omega_{SRF,y}$ on $X_y$ cohomologous to $\omega_0$ where $\omega_0$ is a K\"ahler metric attached to $X$. So there is a smooth function $\rho_y$ on $X_y$
such that $\omega_0\mid_{X_y} +\sqrt{-1}\partial\bar{\partial}\rho_y=\omega_{SRF,y}$ is the unique Ricci-flat K\"ahler metric on
$X_y$. If we normalize $\rho_y$, then
$\rho_y$
varies smoothly in $y$ and defines a smooth function $\rho$ on $X$ and we let

$$\omega_{SRF}=\omega_0+\sqrt{-1}\partial\bar{\partial}\rho$$
which is called as semi-Ricci flat metric. Robert Berman and Y.J.Choi independently showed that the semi-Ricci flat metric is semi positive along horizontal direction. Now for semi Ricci flat metric $\omega_{SRF}$, we have

$$\omega_{SRF}^{n+1}=c(\omega_{SRF}).\omega_{SRF}^ndy\wedge d\bar y$$
Here $c(\omega_{SRF})$ is called a geodesic curvature of semi  $\omega_{SRF}$. Now from Berman and Choi formula, for $V\in T_yY$, the following PDE holds on $X_y$

$$-\Delta_{\omega_{SRF}}c(\omega_{SRF})(V)=|\bar\partial V_{\omega_{SRF}}|_{\omega_{SRF}}^2-\Theta_{V\bar V}(\pi_*(K_{X/Y}))$$

$\Theta_{V\bar V}$ is the Ricci curvature of direct image of relative line bundle( which is a line bundle, since fibers are Calabi Yau manifolds ). Now by integrating on both sides of this PDE, since $$\int_X\Delta_{\omega_{SRF}}c(\omega_{SRF})(V)=0$$
and from the definition of Weil-Petersson metric and this PDE we get $\pi^*\omega_{WP}=Ric(\pi_*(K_{X/Y}))$ and hence for some holomorphic $(n,0)$-form(as non vanishing holomorphic section on the direct image of relative line bundle, which is still line bundle, since fibres are Calabi-Yau varieties)  $\Omega_y$
on $\pi^{-1}(U)$, where $U$ is a neighborhood of $y$ we have 

$$Ric(\pi_*(K_{X/Y}))=-\sqrt[]{-1}\partial\bar\partial \log\|\Omega\|_y^2$$ 

From definition of pushforward of a current by duality, for any continuous function $\psi$ on $Y$, we have

$$\int_{Y}\psi f_*\Omega=\int_{X}(f^*\psi)\Omega=\int_{y\in Y}\int_{\pi^{-1}(y)}(f^*\psi)\Omega$$
and hence on regular part of $Y$ we have

$$\pi_*\Omega=\int_{\pi^{-1}(y)}\Omega $$
and so $$\mid\mid\Omega\mid\mid_y^2=\int_{\pi^{-1}(y)}\mid\Omega\mid_y^2$$
Hence $$\omega_{WP}=-\sqrt[]{-1}\partial_y\bar\partial_y\log\int_{X_y}|\Omega_y|^2$$ and we obtain the desired result.

\;

\;

Now, we give Bismut, Gillet and Christophe Soulé \cite{77},\cite{78} observation about Grothendieck-Riemann-Roch and Chern form of relative canonical bundle which is related to degeneration of K\"ahler-Einstein metric.

One of the applications of Grothendieck-Riemann-Roch is study of degeneration of K\"ahler-Einstein metric. In fact, from Bismut-Gillet-Soulé formula we can get a relation between Weil-Petersson metric on the degeneration of K\"ahler Einstein metric  and fiberwise K\"ahler-Einstein metric by using fiberwise Integral. We explain with more detail.

Let $X$ be a smooth quasi-projective scheme over a field and $K_{0}(X)$ be the Grothendieck group of bounded complexes of coherent sheaves. Consider the Chern character as a functorial transformation
$${\mbox{ch}}\colon K_{0}(X)\to A(X,{\mathbb {Q}})$$, where $A_{d}(X,{\mathbb {Q}})$
is the Chow group of cycles on $X$ of dimension $d$ modulo rational equivalence, tensored with the rational numbers. 
Now consider a proper morphism

   $$f\colon X\to Y$$

between smooth quasi-projective schemes and a bounded complex of sheaves ${{\mathcal {F}}^{\bullet }}$ on $X$. Let $R$ be the right derived functor, then we have the following theorem of Grothendieck-Riemann-Roch.

The Grothendieck-Riemann-Roch theorem relates the pushforward map
    $$f_{\mbox{!}}=\sum (-1)^{i}R^{i}f_{*}\colon K_{0}(X)\to K_{0}(Y)$$

and the pushforward $f_{*}\colon A(X)\to A(Y)$
by the formula
    $${\mbox{ch}}(f_{\mbox{!}}{\mathcal {F}}^{\bullet }){\mbox{td}}(Y)=f_{*}({\mbox{ch}}({\mathcal {F}}^{\bullet }){\mbox{td}}(X))$$
Here $td(X)$ is the Todd genus of (the tangent bundle of) $X$. In fact, since the Todd genus is functorial and multiplicative in exact sequences, we can rewrite the Grothendieck-Riemann-Roch formula as
    $${\mbox{ch}}(f_{\mbox{!}}{\mathcal {F}}^{\bullet })=f_{*}({\mbox{ch}}({\mathcal {F}}^{\bullet }){\mbox{td}}(T_{f}))$$
where $T_f$ is the relative tangent sheaf of $f$, defined as the element $TX - f^*TY$ in $K^0(X)$. For example, when $f$ is a smooth morphism, $T_f$ is simply a vector bundle, known as the tangent bundle along the fibers of $f$.

Let $\pi : \mathcal X \to \mathcal B$ be a proper, smooth holomorphic map of reduced complex spaces. Let $(\mathcal E, \mathcal h)$ be a hermitian vector bundle on
$\mathcal X$ . Let $\lambda(\mathcal E)=\det f_!(\mathcal E)$ be determinant line bundle in the sense of in the derived category.

There exists a Quillen metric $h_Q$ on the determinant line bundle(due to Bismut, Gillet and Soule) such
that the following equality holds for its first Chern form on the base $\mathcal B$ and the component in degree two of a fiber integral:

$$c_1(\lambda(\mathcal E), h^Q)=\left[\int_{\mathcal X/\mathcal B}td(\mathcal X/\mathcal B, \omega_{\mathcal X/\mathcal B})ch(\mathcal E,h)\right]_2$$

Here $ch$ and $td$ stand for the Chern and Todd character forms and $\omega_{\mathcal X/\mathcal B}$ is the fiberwise Kahler form

Now let $(L, h)$ be a hermitian line bundle and take $\mathcal E=(L-L^{-1})^{n+1}$

then the following is known result of  Bismut, Gillet and Soule $$c_1(\lambda(\mathcal E), h^Q)=2^{n+1}\int_{\mathcal X/\mathcal B}c_1(L,h)^{n+1}$$

Now if we take $L=K_{\mathcal X/\mathcal B}=K_{\mathcal X}\otimes \pi^* K_{\mathcal B^{-1}}$ in previous formula, 

If you take fibers are of general type i.e., $c_1({X_b})<0$ , then the left hand side of previous formula is canonical metric on moduli space of fibers which is Weil-Petersson metric. So, we can write

$$\omega_{WP}=-\int_{X/Y}\omega_{SKE}^{n+1}$$

where $\omega_{SKE}$ is the fiberwise K\"ahler-Einstein metric

Moreover if fibers are polarized log Calabi-Yau varieties, then from Bismut-Gillet-Soule formula we have, we can write the logarithmic Weil-Petersson metric as 

$$\omega_{WP}^D=\int_{X'/Y}(\omega_{SRF}^D)^{n+1}$$

where $\omega_{SRF}^D$ is the fiberwise Ricci flat metric( where $X'=X\setminus D$).

or fibers be K-stable then we have the same formula for canonical metric on moduli space of K-stable fibers which we still call them Weil-Petersson metric.

$$\omega_{WP}=\int_{X/Y}\omega_{SKE}^{n+1}$$

where $\omega_{SKE}$ is the fiberwise K\"ahler-Einstein metric on moduli space of K-stable Fano fibers which admit K\"ahler-Einstein metric with positive Ricci curvature. See \cite{83}, \cite{79}.

\section*{Logarithmic Weil-Petersson metric via Canonical bundle formula}

Let $f:X\to B$ be the surjective morphism of a normal projective variety $X$
of dimension $n=m+l$ to a nonsingular projective $l$-fold $B$ such that:

i) $(X,D)$ is sub klt pair.

ii) the generic fiber
$F$ of $\pi$ is a geometrically irreducible variety with vanishing log Kodaira dimension. We fix the smallest
$b\in \mathbb Z>0$ such that the

$$\pi_*\mathcal O_X(b(K_X+D))\neq 0$$

The Fujino-Mori \cite{25} log-canonical bundle formula
for $\pi:(X,D)\to B$ is 

$$b(K_X+D)=\pi^*(bK_B+L_{(X,D)/B}^{ss})+\sum_Ps_P^D\pi^*(P)+B^D$$
where $B^D$ is $\mathbb Q$-divisor on $X$ such that $\pi_*\mathcal O_X([iB_+^D])=\mathcal O_B$ ($\forall i>0$). Here $s_P^D:=b(1-t_P^D)$ where $t_P^D$ is the log-canonical threshold of $\pi^*P$ with respect to $(X,D-B^D/b)$ over the generic point $\eta_P$ of $P$. i.e., 

$$t_P^D:=\max \{t\in \mathbb R\mid \left(X,D-B^D/b+t\pi^*(P)\right)\;  \text{is sub log canonical over}\; \eta_P\}$$

 Let $\pi:(X,D)\to B$
is a holomorphic submersion onto
a compact K\"ahler manifold
$B$ with $B$ be a Calabi-Yau manifold, log fibers $(X_s,D_s)$ are log Calabi-Yau, and $D$ is a simple normal crossing divisor in $X$ with conic singularities.
Then $(X,D)$
admits a singular hermitian 
metric $h$ on direct image of relative line bundle,
solving 

$$Ric(h)=\omega_{WP}^D+\sum_P(b(1-t_P^D))[\pi^*(P)]+[B^D]$$
as current where the singular hermitian metric $h$ has zero Lelong number and is good metric in the sense of Mumford. So the Canonical bundle formula gives a nice picture of logarithmic Weil-Petersson metric on moduli space of log Calabi-Yau pairs.

Now we give a motivation that why the geometry of pair $(X,D)$ must be interesting. The first one comes from algebraic geometry, in fact for deforming the cone angle we need to use of geometry of pair $(X,D)$. In the case of minimal general type manifold the canonical bundle of $X$, i.e., $K_X$ is nef and we would like $K_X$ to be ample and it is not possible in general and what we can do is that to add a small multiple of ample bundle $\frac{1}{m}A$, i.e., $K_X+\frac{1}{m}A$ and then we are deal with the pair $(X,\frac{1}{m}H)$ which $H$ is a generic section of it. The second one is the works of Chen-Sun-Donaldson and Tian on existence of K\"ahler Einstein metrics for Fano varieties which they used of geometry of pair $(X,D)$ for their proof .\cite{44}

Now we explain {Tian-Yau program} to how to construct model metrics in general, like conical model metric, Poincare model metric, or Saper model metric,due to A.Fujiki's method

Let $\mathbb C^n=\mathbb C^n(z_1,...,z_n)$ be a complex Euclidian space for some $n > 0$. 
For a positive number $\epsilon$ with $0 < \epsilon < 1$ consider $$\overline X=\overline X_\epsilon=\{z=(z_1,...,z_n)\in \mathbb C^n|\;\; |z_i|<\epsilon\}$$

Now, let $D_i=\{z_i=0\}$ be the irreducible divisors and take $D=\sum_iD_i$ where $$D=\{z\in \overline X| \;\; z_1z_2...z_k=0\}$$ and take $X=\overline X\setminus D$. In polar coordinate we can write $z_i=r_ie^{i\theta_i}$.Let $g$ be a K\"ahler metric on $D$ such that the associated K\"ahler form $\omega$ is of the following form

$$\omega=\sqrt[]{-1}\sum_i\frac{1}{|dz_i|^2}dz_i\wedge d\bar z_i$$

Then the volume form $dv$ associated to $\omega$ is written in the form;

$$dv=(\sqrt[]{-1})^n\prod_{i=1}^n\frac{1}{|dz_i|^2}\prod_idz_i\wedge d\bar{z_i}\;\;,  \; v=\frac{1}{|dz_i|^2}$$

Let $L$ be a (trivial) holomorphic line bundle defined on $\overline X$, with a generating holomorphic section $S$ on $\overline X$. Fix a $C^\infty$ hermitian metric $h$ of $L$ over $X$ and denote by
$|S|^2$ the square norm of $S$ with respect to $h$. Assume the functions $|S|^2$ and $|dz_i|^2$ depend only on $r_i$,\; $1\leq i\leq k$. Set $$d(r_1,...,r_k)=|S|^2.v.\prod_{1\leq i\leq k}r_i$$

and further make the following three assumptions:

\textbf{A1}) The function $d$ is of the form

$$d(r_1,...,r_k)=r_1^{c_1}...r_k^{c_k}(\log 1/r_1)^{b_1}....(\log 1/r_k)^{b_k}L(r_1,...,r_k)^t$$

where $$L=L(r_1,...,r_k)=\sum_{i=1}^k\log 1/r_i$$

and $c_i, b_j, t$ are real numbers with $t \geq 0$ such that $q_i=b_i+t\neq -1$ if $c_i$ is
an odd integer. We set $a_i= (c_i + 1)/2$ and denote by $[a_i]$ the largest integer which does not exceed $a_i$.

\textbf{A2}) If $1\leq i\leq k$, then $|dz_i|^2$ is either of the following two forms;

$$|dz_i|^2(r)=r_i^2(\log 1/r_i)^2, \;\;\text{or}\;\; |dz_i|^2(r)=r_i^2L^2$$

In fact, A2) implies that the K\"ahler metric $g$ is (uniformly) complete along $D$.

\textbf{A3}) If  $k+ 1 \leq i \leq n$, then $|dz_i|^{-2}$ is bounded (above) on $X$.

Now, we give some well-known examples of Tian-Yau-Fujiki picture, i.e, conical model metric, Poincare model metric, and Saper model metric.

 A K\"ahler current $\omega$
is called a conical K\"ahler metric (or Hilbert Modular type) with angle $2\pi\beta$, $(0< \beta <1)$ along the divisor $ D$, if
$\omega$ is smooth away from $D$ and
asymptotically equivalent along $D$ to the model conic metric 

$$\omega_{\beta}=\sqrt{-1}\left(\frac{dz_1\wedge d\bar{z_1}}{|z_1|^{2(1-\beta)}}+\sum_{i=2}^n dz_i\wedge d\bar{z_i}\right)$$
 here $(z_1,z_2,...,z_n)$ are local holomorphic coordinates and $D=\{z_1=0\}$ locally.

After an appropriate -singular- change of coordinates, one can see that this
model metric represents an Euclidean cone of total angle $\theta=2\pi\beta$, whose
model on $\mathbb R^2$ is the following metric: $d\theta^2 +\beta^2dr^2$.
The volume form $V$ of a conical K\"ahler metric $\omega_D$ on the pair $(X,D)$ has the form
$$V=\prod_j|S_j|^{2\beta_j-2}e^f\omega^n$$ where $f\in C^0$.

This asymptotic behaviour of metrics can be translated to the second order
asymptotic behaviour of their potentials

$$\omega_{\beta}=-\sqrt{-1}\partial\bar{\partial}\log e^{-u}$$

where $u=\frac{1}{2}\left(\frac{1}{\beta^2}|z_1|^{2\beta}+|z_2|^2+...+|z_n|^2\right)$.

Moreover, if we let $z=re^{i\theta}$ and $\rho=r^\beta$ then the model metric in $\omega_\beta$ becomes

$$(d\rho+\sqrt{-1}\beta\rho d\theta)\wedge (d\rho-\sqrt{-1}\beta\rho d\theta)+\sum_{i>1}dz_i\wedge d\bar{z_i}$$

and if we set $\epsilon=e^{\sqrt{-1}\beta\theta}(d\rho+\sqrt{-1}\beta\rho d\theta)$ then
the conical K\"ahler metric $\omega$ on $(X,(1-\beta)D)$ can be expressed as 

$$\omega=\sqrt{-1}\left(f\epsilon\wedge\bar\epsilon+f_{\bar j}\epsilon\wedge d\bar z_j+f_jdz_j\wedge\bar\epsilon+f_{i\bar j}dz_i\wedge d\bar z_j\right)$$

By the assumption on the asymptotic behaviour we we mean there exists some coordinate chart in which the zero-th order asymptotic of the metric agrees with the model metric. In other words,
there is a constant $C$, such that

$$\frac{1}{C}\omega_\beta\leq\omega\leq C\omega_\beta$$

\:
\:
\:

In this note because we assume certain singularities for the K\"ahler manifold $X$ we must design our K\"ahler Ricci flow such that our flow preserve singularities. Now fix a simple normal crossing divisor $D=\sum_{i} (1-\beta_i)D_i$, where $\beta_i\in (0,1)$ and simple normal crossing divisor $D$ means that $D_i$'s are irreducible smooth divisors and for any $p\in Supp (D)$ lying in the intersection of exactly $k$ divisors $D_1,D_2,...,D_k$ , there exists a coordinate chart $(U_p, \{z_i\})$ containing $p$, such that $D_i|_{U_p}=\{z_i=0\}$ for $i=1,...,k$.

If  $S_i\in H^0(X,\mathcal O_X\left(L_{D_i}\right))$ is the defining sections and $h_i$ is  hermitian  metrics on the line bundle induced
by $D_i$, then Donaldson showed that for sufficiently small $\epsilon_i>0$, $\omega_i=\omega_0+\epsilon_i\sqrt{-1}\partial\overline{\partial}|S_i|_{h_i}^{2\beta_i}$ gives a conic K\"ahler metric on $X\setminus Supp(D_i)$ with cone angle $2\pi \beta_i$ along divisor $D_i$ and also if we set $\omega=\sum_{i=1}^N\omega_i$ then,
$\omega$ is a smooth K\"ahler metric on $X\setminus Supp(D)$ and $$||S||^{2(1-\beta)}=\prod_{i=1}^k||S_i||^{2(1-\beta)}$$ where $S\in H^0(X,\mathcal O(L_D))$. Moreover, $\omega$ is uniformly equivalent to the standard cone metric 

$$\omega_p=\sum_{i=1}^k\frac{\sqrt{-1}dz_i\wedge d\bar{z_j}}{|z_i|^{2(1-\beta_i)}}+\sum_{i=k+1}^N\sqrt{-1}dz_i\wedge d\bar{z_i}$$

From Fujiki theory, $|dz_i|^2=r_i^2$ for $1\leq i\leq k$ and $|dz_j|^2=1$ for $k+1\leq j\leq n$ so that
A2) and A3) are again satisfied.

From now on for simplicity we write just "divisor $D$" instead "simple normal crossing divisor $D$".

We give an example of varieties which have conical singularities.  Consider a smooth geometric orbifold given by $\mathbb Q$-divisor $$D=\sum_{j\in J}(1-\frac{1}{m_j})D_j$$ where $m_j\geq 2$ are positive integers and $\text{Supp}D=\cap_{j\in J}D_j$ is of normal crossings divisor. Let $\omega$ be any K\"ahler metric on $X$, let
$C >0$ be a real number and $s_j\in H^0\left(X,\mathcal O_X(D_j)\right)$ be a section defining $D_j$. Consider the following expression 

$$\omega_D=C\omega+\sqrt[]{-1}\sum_{j\in J}\partial\bar\partial |s_j|^{2/m_j}$$

If $C$ is large enough, the above formula defines a closed positive $(1,1)$
-current (smooth away from $D$). Moreover

$$\omega_D\geq \omega$$ in the sense of currents. Consider $\mathbb C^n$ with the orbifold divisor given by the equation

$$\prod_{j=1}^nz_j^{1-{1}/{m_j}}=0$$ (with eventually $m_j=1$ for some $j$). The sections
$s_j$ are simply the coordinates $z_j$ and a simple computation gives

$$\omega_D=\omega_{eucl}+\sqrt[]{-1}\sum_{j=1}^n\partial\bar\partial |z_j|^{2/m_j}=\omega_{eucl}+\sqrt[]{-1}\sum_{j=1}^n\frac{dz_j\wedge d\bar {z_j}}{m_j^2|z_j|^{2(1-1/m_j)}}$$

Here we mention also metrics with non-conic singularities. We say a metric $\omega$ is of Poincare type, if it is quasi-isometric to 
$$\omega_\beta=\sqrt{-1}\left(\sum_{i=1}^k\frac{dz_i\wedge d\bar z_i}{|z_i|^2\log^2|z_i|^2}+\sum_{i=k+1}^n dz_i\wedge d\bar z_i\right)$$
It is always possible to construct a Poincare metric on $ M \setminus D$ by
patching together local forms with $C^{\infty}$ partitions of unity. Now, from Fujiki theory $|dz_i|^2=r_i^2(\log 1/r_i)^2$, $1\leq i\leq k$ and $|dz_j|^2=1$, $k+1\leq j\leq n$ so that A2) and A3) above are satisfied; we have

$$v=\prod_{1\leq i\leq k}r_i^{-2}(\log 1/r_i)^{-2}$$

\;

Let $\Omega_P$ be the volume form on $X\setminus D$, then, there exists a locally bounded positive continuous function $c(z)$ on polydisk $\mathbb D^n$ such that 

$$\Omega_P=c(z)\sqrt{-1}\left(\wedge_{i=1}^k\frac{dz_i\wedge d\bar z_i}{|z_i|^2\log^2|z_i|^2}+\wedge_{i=k+1}^n dz_i\wedge d\bar z_i\right)$$

holds on $\mathbb D^n\cap (X\setminus D)$

\textbf{Remark A}: Note that if $\Omega_P$ be a volume form of Poincare growth on $(X,D)$, with $X$ compact. If $c(z)$ be $C^2$ on $\mathbb D^n$, then $-Ric(\Omega_P)$ is of Poincare growth.

We say that $\omega$ is the homogeneous Poincare metric if its fundamental form $\omega_\beta$ is described
locally in normal coordinates by the quasi-isometry

$$\omega_\beta=\sqrt{-1}\left(\frac{1}{(\log|z_1z_2...z_k|^2)^2}\sum_{i=1}^k\frac{dz_i\wedge d\bar z_i }{|z_i|^2}+\sum_{i=1}^ndz_i\wedge d\bar z_i\right)$$
and we say $\omega$ has Ball Quotient singularities if it is quasi-isometric to 

$$\omega_\beta=\sqrt{-1}\frac{dz_1\wedge d\bar{z_1}}{(|z_1|\log (1/|z_1|))^2}+\sqrt{-1}\sum_{j=2}^n\frac{dz_j\wedge d\bar z_j}{\log 1/|z_1|}$$
It is called also Saper's distinguished metrics.

$$|dz_1|^2=r_1^2(\log 1/r_1)^2,\;\;\; |dz_j|^2=\log 1/r_1,\;\; k+1\leq j \leq n$$ so that A2) and
A3) are satisfied; also we have the volume form as

$$v=r_1^{-2}\left(\log 1/r_1\right)^{-(n+1)}$$

If $\omega$ is the fundamental form of a metric on the compact manifold $X$, and $\omega_{sap}$ be the fundamental forms of Saper's distinguished metrics and $\omega_{P,hom}$ be the fundamental forms of homogeneous Poincare metric, 
 on the noncompact manifold $M\setminus D$, then $\omega_{sap}+\omega$ and $\omega_{P,hom}$ are quasi-isometric.

For the log-Calabi-Yau fibration $f:(X,D)\to Y$, such that $(X_t,D_t)$ are log Calabi-Yau varieties and central fiber $(X_0,D_0)$ has simple normal crossing singularities,and if $(X,\omega)$ be a K\"ahler variety with Poincaré singularities then the semi Ricci flat metric $\omega_{SRF}$ has pole singularities with Poincare growth. So the semi-Ricci flat metric $\omega_{SRF}|_{X_t}$ is quasi-isometric with the following model which we call it fibrewise Poincaré singularities.

$$\frac{\sqrt[]{-1}}{\pi}\sum_{k=1}^n\frac{dz_k\wedge d\bar {z_k}}{|z_k|^2(\log|z_k|^2)^2}+\frac{\sqrt[]{-1}}{\pi}\frac{1}{\left(\log|t|^2-\sum_{k=1}^n\log|z_k|^2\right)^2}\left(\sum_{k=1}^n\frac{dz_k}{z_k}\wedge \sum_{k=1}^n\frac{d\bar {z_k}}{\bar {z_k}}\right)$$

We can define the same fibrewise conical singularities and the semi-Ricci flat metric $\omega_{SRF}|_{X_t}$ is quasi-isometric with the following model

$$\frac{\sqrt[]{-1}}{\pi}\sum_{k=1}^n\frac{dz_k\wedge d\bar {z_k}}{|z_k|^2}+\frac{\sqrt[]{-1}}{\pi}\frac{1}{\left(\log|t|^2-\sum_{k=1}^n\log|z_k|^2\right)^2}\left(\sum_{k=1}^n\frac{dz_k}{z_k}\wedge \sum_{k=1}^n\frac{d\bar {z_k}}{\bar {z_k}}\right)$$

Note that if the base of fibration be smooth and of dimension one, then always by using base change and birational change we can assume the central fiber has simple normal crossing singularities due to following Mumford et al. theorem. Note that if we assume the general fibers are log-Calabi-Yau pairs, then we can assume the central fiber is Calabi-Yau pair after using semi-stable minimal model program \cite{69}.

\;

\textbf{Theorem} (Grothendieck, Kempf, Knudsen, Mumford and Saint-Donat\cite{70})Let $ k$ be an algebraically closed field of characteristic 0 (e.g.  $ k={\mathbb{C}}$). Let $ f:X\to C$ be a surjective morphism from a  $ k$-variety $ X$ to a non-singular curve $ C$ and assume there exists a closed point $ z\in C$ such that  $ f_{\vert X\setminus f^{-1}(z)}:
X\setminus f^{-1}(z)\to C\setminus\{z\}$ is smooth. Then we find a commutative diagram

$$\displaystyle \xymatrix{ X\ar[d]_f & X\times_C C'\ar[l]\ar[d] &
X'\ar[l]_-{p}\ar[dl]^{f'}\\
C & C'\ar[l]^{\pi}
}$$
with the following properties

\;

1. $ \pi:C'\to C$ is a finite map, $ C'$ is a non-singular curve and  $ \pi^{-1}(z)=\{z'\}$.

\;

2. $ p$ is projective and is an isomorphism over  $ C'\setminus
\{z'\}$.
$ X'$ is non-singular and  $ {f'}^{-1}(z')$ is a reduced divisor with simple normal crossings, i.e., we can write  $ {f'}^{-1}(z')=\sum_i E_i$ where the $ E_i$ are 1-codimensional subvarieties (i.e., locally they are defined by the vanishing of a single equation), which are smooth and, for all  $ r$, all the intersections  $ E_{i_1}\cap\ldots\cap E_{i_r}$ are smooth and have codimension $ r$.

\textbf{Remark}: For a polarized family of varieties over $\mathbb D$, smooth over
$\mathbb D^\times$ and with smooth fibers Calabi-Yau, Tosatti and Takayama \cite{80}, \cite{81} showed that the origin lies at
finite Weil-Petersson distance if and only if after a finite base change the
family is birational to one with central fiber a Calabi-Yau variety with at worst canonical singularities. It is worth to mention that we can have the same result when we replace the finiteness of Weil-Petersson distance with boundedness of Song-Tian-Tsuji measure, i.e, Song-Tian-Tsuji measure is bounded near origin if and only if after a finite base change the
family is birational to one with central fiber a Calabi-Yau variety with at worst canonical singularities. In fact if we consider the degeneration of Calabi-Yau varieties when the dimension of base is bigger than one, we need to replace semi-stable reduction of Mumford with weak semi-stable Abramovich-Karu reduction\cite{82} . So the Song-Tian-Tsuji measure is more suitable instead of Weil-Petersson distance.

\textbf{Corollary}: Wang \cite{76} showed that , If $X$ be a Calabi-Yau varieties which admits a smoothing to Calabi-Yau manifolds. If $X$ has only canonical singularities then $X$ has finite
Weil-Petersson distance along any such smoothing. Now we can extend this result and in fact  If $X$ be a Calabi-Yau varieties which admits a smoothing to Calabi-Yau manifolds and $X$ has only canonical singularities then Song-Tian-Tsuji measure is bounded along any such smoothing.

\begin{defn}
A K\"ahler metric with cone singularities along
$D$ with cone angle $2\pi\beta$
is a smooth
K\"ahler metric on
$X\setminus D$
which satisfies the following conditions when we write $\omega_{sing}=\sum_{i,j}g_{i\bar j}\sqrt{-1}dz_i\wedge d\bar z_j$ in terms of the local holomorphic coordinates $(
z_1;...; z_n)$ on a neighbourhood $U\subset X$
with $D\cap U=\{z_1=0\} $

\;

1. $g_{1\bar 1}=F|z_1|^{2\beta-2}$ for some strictly positive smooth bounded function
$F$
on
$X
\setminus
D$

\;

2. $g_{1\bar j}=g_{i\bar 1}=O(|z_1|^{2\beta-1})$

\;

3. $g_{i\bar j}=O(1)$ for
$i,j
\neq 1$

\end{defn}

Now we shortly explain Donaldson's linear theory which is useful later in the definition of semi ricci flat metrics.

\begin{defn}
1) A function $f$ is in $C^{,\gamma,\beta} (X, D)$ if $f$ is $C^{\gamma}$ on $X \setminus D$, and locally near
each point in $D$, $f$ is $C^{\gamma}$ in the coordinate 
$( \hat\zeta =\rho e^{i\theta} = z_1|z_1|^{\beta-1}, z_j )$.

\;

2)A (1,0)-form $\alpha$ is in $C^{,\gamma,\beta} (X, D)$ if $\alpha$ is $C^\gamma$ on $X \setminus D$ and locally near each point in $D$, we have $\alpha = f_1\epsilon +\sum_{ j>1} f_jdz_j$ with $f_i \in C^{,\gamma,\beta}$ for $1 \leq i \leq n$, and $f_1 \to 0$ as $z_1 \to 0$ where $\epsilon=e^{\sqrt{-1}\beta\theta}(d\rho+\sqrt{-1}\beta\rho d\theta)$

\;

3) A $(1,1)$-form $\omega$ is in $C^{,\gamma,\beta} (X, D)$ if $\omega$ is $C^\gamma$ on $X \setminus D$ and near each point in
$D$ we can write $\omega$ as $$\omega=\sqrt{-1}\left(f\epsilon\wedge\bar\epsilon+f_{\bar j}\epsilon\wedge d\bar z_j+f_jdz_j\wedge\bar\epsilon+f_{i\bar j}dz_i\wedge d\bar z_j\right)$$

 such that
$f,f_j,f_{\bar j},f_{i\bar j}\in C^{,\gamma,\beta}$, and $f_j,f_{\bar j}\to 0$ as $z_1\to 0$

\;

4)A function $f$ is in $C^{2,\gamma,\beta} (X, D)$ if $f$,$\partial f$,$\partial\bar\partial f$
are all in
$C^{,\gamma,\beta} $

\end{defn}

Fix a smooth metric $\omega_0$ in $c_1(X)$, we define the space of admissible functions to be

$$\hat C(X,D)=C^{2,\gamma}(X)\cup
\bigcup_{0<\beta<1} \left ( \bigcup_{0<\gamma<\beta^{-1}-1}C^{2,\gamma,\beta}(X,D)\right)$$
and the space of admissible K\"ahler potentials to be

$$\hat {\mathcal H }(\omega_0)=\{\phi\in \hat C(X,D)\mid \omega_\phi=\omega_0+\sqrt{-1}\partial\bar\partial \phi>0\}$$

Note that

$$ {\mathcal H }(\omega_0)\subset \hat {\mathcal H }(\omega_0)\subset \mathcal P\mathcal S{\mathcal H }(\omega_0)\cap L^{\infty}(X)$$

Where $\mathcal P\mathcal S{\mathcal H }(\omega_0)\cap L^{\infty}(X)$ is the space of bounded $\omega_0$-plurisubharmonic functions and 

$$\mathcal P\mathcal S{\mathcal H }(\omega_0)=\{\phi\in L_{loc}^1(X)\mid \phi \;  \text{is u.s.c and}\; \omega_0+\sqrt{-1}\partial\bar\partial \phi>0\}$$

The Ricci curvature of the K\"ahlerian form $\omega_D$ on the pair $(X,D)$ can be represented as:

$$Ric\left(\omega_D\right)=2\pi\sum_j(1-\beta_j)[D_j]+\theta+\sqrt{-1}\partial\bar{\partial}\psi$$
with $\psi\in C^0(X)$ and $\theta$ is closed smooth $(1,1)$-form.

We have also $dd^c$-lemma on $X=\overline X\setminus D$. Let $\Omega$ be a smooth closed $(1,1)$-form in the cohomology class $c_1(K_{\overline X}^{-1}\otimes L_D^{-1})$. Then for any $\epsilon>0$ there exists an explicitly given complete K\"ahler metric $g_\epsilon$ on $M$ such that 

$$Ric(g_\epsilon)-\Omega=\frac{\sqrt[]{-1}}{2\pi}\partial\bar\partial f_\epsilon\;\;  \text{on} X$$
where $f_\epsilon$ is a smooth function on $X$ that decays to the order of $O(\|S\|^\epsilon)$. Moreover, the Riemann curvature tensor $R(g_\epsilon)$ of the metric $g_\epsilon$ decays to the order of $O\left((-n\log \|S\|^2)^{-\frac{1}{n}}\right)$

 Now we explain the logarithmic Weil-Petersson metric on moduli space of log Calabi-Yau manifolds(if it exists. for special case of rational surfaces it has been proven that such moduli space exists). The logarithmic Weil-Petersson metric has pole singularities and we can introduce it also by elements of logarithmic Kodaira-Spencer tensors which represent elements of $H^1\left(X,\Omega_X^1(\log(D))^\vee\right)$. More precisely, Let $X$ be a complex manifold, and $D \subset X$ a divisor and $\omega$ a holomorphic $p$-form on $X\setminus D$. If $\omega$ and $d\omega$ have a pole of order at most one along $D$, then $\omega$ is said to have a logarithmic pole along $D$. $\omega$ is also known as a logarithmic $p$-form. The logarithmic $p$-forms make up a subsheaf of the meromorphic $p$-forms on $X$ with a pole along $D$, denoted

$$\Omega^p_X(\log D)$$

and for the simple normal crossing divisor $D=\{z_1z_2...z_k=0\}$ we can write the stalk of $\Omega^1_X(\log D)$ at $p$ as follows

$$\Omega_X^1(\log D)_p = \mathcal{O}_{X,p}\frac{dz_1}{z_1}\oplus\cdots\oplus\mathcal{O}_{X,p}\frac{dz_k}{z_k} \oplus \mathcal{O}_{X,p}dz_{k+1} \oplus \cdots \oplus \mathcal{O}_{X,p}dz_n$$

Since, fibers are log Calabi-Yau manifolds and by recent result of Jeffres-Mazzeo-Rubinstein, we have Ricci flat metric on each fiber $(X_y,D_y)$ and hence we can have log semi-Ricci flat metric and by the same method of previous theorem, the proof of Theorem 8 is straightforward.

\begin{thm}
Let $(M,\omega_0)$ be a compact K\"ahler manifold with $D\subset M$ a smooth divisor and suppose we have topological constraint condition $c_1(M)=(1-\beta)[D]$ where $\beta\in (0,1]$ then there exists a conical K\"ahler Ricci flat metric with angle $2\pi\beta$ along $D$. This metric is unique in its K\"ahler class. This metric is polyhomogeneous; namely, the K\"ahler Ricci flat metric $\omega_0+\sqrt{-1}\partial\bar\partial \varphi$ admits a complete asymptotic expansion with smooth coefficients as $r\to 0$ of the form 

$$\varphi(r,\theta,Z)\sim \sum_{j,k\geq 0}\sum_{l=0}^{N_{j,k}}a_{j,k,l}(\theta,Z)r^{j+k/\beta}(\log r)^l$$

where $r=|z_1|^\beta/\beta$ and $\theta=\arg z_1$ and with each $a_{j,k,l}\in C^\infty$. 

\end{thm}

Now we can introduce Logarithmic semi Ricci flat metrics. The volume of fibers $(X_y,D_y)$ are homological constant independent of $y$, and we assume that it is equal to $1$. Since fibers are log Calabi-Yau varieties, so $c_1(X_y,D_y)=0$, hence there is a smooth function $F_y$ such that $Ric(\omega_y)=\sqrt{-1}\partial\bar{\partial}F_y$ . The function $F_y$ vary smoothly in $y$. By Jeffres-Mazzeo-Rubinstein's theorem \cite{61}, there is a unique conical Ricci-flat K\"ahler metric $\omega_{SRF,y}$ on $X_y\setminus D_y$ cohomologous to $\omega_0$. So there is a smooth function $\rho_y$ on $X_y\setminus D_y$
such that $\omega_0\mid_{X_y\setminus D_y} +\sqrt{-1}\partial\bar{\partial}\rho_y=\omega_{SRF,y}$ is the unique Ricci-flat K\"ahler metric on
$X_y\setminus D_y$. If we normalize $\rho_y$, then
$\rho_y$
varies smoothly in $y$ and defines a smooth function $\rho^D$ on $X\setminus D$ and we let

$$\omega_{SRF}^D=\omega_0+\sqrt{-1}\partial\bar{\partial}\rho^D$$
which is called as Log Semi-Ricci Flat metric.

Now we explain one of applications of semi-Ricci flat metric in Mirror symmetry in the context of Strominger-Yau and Zaslow conjecture..
In 2002, Hausel- Thaddeus\cite{84}  interpreted SYZ conjecture in the context of Hitchin system and Langlands duality.

Let $\pi : E \to \Sigma$  a complex vector bundle of rank $r$ and degree $d$ equipped with a hermitian metric on Riemann surface $\Sigma$ . Take th moduli space 
$$M(r, d) = \{(A, \Phi) \text{ solving }(\star)\}/\mathcal G $$

(which is a finite-dimensional non-compact space carrying a natural hyper-Kähler  metric)

where 

$$F^0_A + [\Phi \wedge \Phi^*] = 0 ,\; \;  \bar \partial A\Phi =0\; \; (\star)$$

Here $A$ is a unitary connection on $E$ and $\Phi \in \Omega^{1,0}(End E)$ is a Higgs field. $F^0$ denotes the trace-free part of the curvature and $\mathcal G$ is the unitary gauge group.

$M(r, d)$ is the total space of an integrable system(which can be interpreted by the non-abelian Hodge theory due to Corlette), the Hitchin fibration, together with Langlands
duality between Lie groups provides a model for mirror symmetry in the Strominger-Yau and Zaslow conjecture.

In fact rank 2 Hitchin fibration $$\det: M(2, d)\to H^0(\Sigma, K^2_\Sigma),\;  [(A, \Phi)]\to \det \Phi$$ gives rise to semi Ricci flat metric $\omega_{SRF}$ on $ M(2, d)_{reg}$ which still is open to be semi-positive as current. Study of analysis of such semi-Ricci flat metric led to solution of a lot of conjectures in mirror symmetry

Let $f:X\setminus D\to S$, be a smooth family of quasi-projective K\"ahler manifolds. Let $x\in X\setminus D$, and $(\sigma,z_2,...,z_n,s^1,...,s^d)$, be a coordinate centered at $x$, where $(\sigma,z_2,...,z_n)$ is a local coordinate of a fixed fiber of $f$ and $(s^1,...,s^d)$ is a local coordinate of $S$, such that $$f(\sigma,z_2,...,z_n,s^1,...,s^d)=(s^1,...,s^d)$$

Now consider a smooth form $\omega$ on $X\setminus D$, whose restriction to any fiber of $f$, is positive definite. Then $\omega$ can be written as 

\begin{align*}
\omega(\sigma,z,s)&=\sqrt[]{-1}(\omega_{i\bar j}ds^i\wedge d\bar s^{ j}+\omega_{i\bar \beta}ds^i\wedge d\bar z^{ \beta}+\omega_{\alpha\bar j}dz^\alpha\wedge d\bar s^{ j}+\omega_{\alpha\bar \beta}dz^\alpha\wedge d\bar z^{ \beta}+\omega_{\sigma}d\sigma\wedge d\bar s^{ j}\\&+
\omega_{i\bar \sigma}ds^i\wedge d\bar \sigma+\omega_{\sigma\bar\sigma}d\sigma\wedge d\bar \sigma+\omega_{\sigma\bar j}d\sigma\wedge d\bar z^{ j}+\omega_{i\bar \sigma}dz^i\wedge d\bar \sigma)
\end{align*}

Since $\omega$ is positive definite on each fibre, hence 

\begin{align*}\sum_{\alpha,\beta=2}\omega_{\alpha\bar \beta}dz^\alpha\wedge d\bar z^{ \beta}+\omega_{\sigma\bar\sigma}d\sigma\wedge d\bar \sigma+\sum_{j=2}\omega_{\sigma\bar j}d\sigma\wedge d\bar z^{ j}+\sum_{i=2}\omega_{i\bar \sigma}dz^i\wedge d\bar \sigma\end{align*}

gives a K\"ahler metric on each fiber $X_s\setminus D_s$. So

$$\det(\omega_{\lambda\bar\eta}^{-1}(\sigma,z,s))=\det\bordermatrix{&\omega_{\sigma\bar\sigma} & \omega_{\sigma\bar 2}  & \ldots & \omega_{\sigma\bar n}\cr & \omega_{2\bar\sigma}  &  \omega_{2\bar 2} & \ldots & \omega_{2\bar n}\cr
 & \vdots & \vdots & \ddots & \vdots\cr
 & \omega_{n\bar\sigma}  &   \omega_{2\bar n}       &\ldots & \omega_{n\bar n}}^{-1}$$
gives a hermitian metric on the relative line bundle $K_{X'/S}$ and its Ricci curvature can be written as $\sqrt[]{-1}\partial\bar\partial \log\det \omega_{\lambda\bar\eta}(\sigma,z,s)$

\begin{thm}By the same method we can introduce the logarithmic Weil-Petersson metric on $\pi:(X,D)\to Y$ with assuming fibers to be log Calabi-Yau manifolds and snc divisor $D$ has conic singularities, then we have 
$$\omega_{WP}^D=-\sqrt[]{-1}\partial_y\bar\partial_y\log\int_{X_y\setminus D_y}\frac{\Omega_y\wedge\bar\Omega_y}{\|S_y\|^2}$$ 
where $S_y\in H^0(X_y,L_{D_y})$. Moreover, if $\omega$ has Poincare singularities along snc divisor $D$, we have the following formula for logarithmic Weil-Petersson metric

$$\omega_{WP}^D=-\sqrt[]{-1}\partial_y\bar\partial_y\log\int_{X_y\setminus D_y}\frac{\Omega_y\wedge\bar\Omega_y}{\|S_y\|^2\log^2\|S_y\|^2}$$ 

\end{thm}

\begin{thm} Let $\pi:(X,D)\to Y$ be a holomorphic family of log Calabi-Yau pairs $(X_s,D_s)$ for the K\"ahler varieties $X,Y$. Then we have the following relation between logarithmic Weil-Petersson metric and fiberwise Ricci flat metric.

$$\sqrt[]{-1}\partial\bar{\partial}\log (\frac{f^*\omega_Y^m\wedge (\omega_{SRF}^D)^{n-m}}{\mid S\mid^2})=-f^*Ric(\omega_Y)+f^*\omega_{WP}^D$$ where $S\in H^0(X,\mathcal O(L_N))$, here $N$ is a divisor which come from Fujino-Mori's canonical bundle formula(see \cite{38} also)
\end{thm}

\textbf{Proof}: Take $X'=X\setminus D$. Choose a local nonvanishing holomorphic section $\Psi_y$ of $\pi_*(K_{X'/Y}^{\otimes l})$ with $y\in U\subset X'$. We define a smooth positive function on $\pi(U)$ by

$$u(y)=\frac{(\sqrt[]{-1})^{(n-m)^2}(\Psi_y\wedge \overline{\Psi_y})^{\frac{1}{l}}}{(\omega_{SRF}^D)^{n-m}\mid_{X_y\setminus D_y}}$$

But the Numerator and Denominator of $u$ are Ricci flat volume forms on $X_y\setminus D_y$, so $u$ is a constant function. Hence by integrating $u(y)(\omega_{SRF}^D)^{n-m}\mid_{X_y\setminus D_y}$ over $X_y\setminus D_y$ we see that 

$$u(y)=\frac{(\sqrt[]{-1})^{(n-m)^2}\int_{X_y\setminus D_y}\frac{(\Psi_y\wedge \overline{\Psi_y}}{\mid S_y\mid^2})^{\frac{1}{l}}}{\int_{{X_y\setminus D_y}}\frac{(\omega_{SRF}^D)^{n-m}}{\mid S_y\mid^2}}$$

where $S_y\in H^0(X',\mathcal O(L_{D_y}))$.

But $y\mapsto \int_{X_y\setminus D_y}\frac{(\omega_{SRF}^D)^{n-m}}{\mid S_y\mid^2}$ is constant over $Y$. Hence the Logarithmic Weil-Petersson can be written as $$-\sqrt[]{-1}\partial\bar{\partial}\log u=\omega_{WP}^D\;\;\;(*)$$

Now, to finish the proof we can write $\Psi
_y=F(\sigma,y,z)(d\sigma\wedge dz_2\wedge...\wedge dz_{n-m})$ where $F$ is holomorphic and non-zero. Hence by substituting $\Psi
_y$ in $u$ and rewriting $\sqrt[]{-1}\partial\bar{\partial}\log (\frac{f^*\omega_Y^m\wedge (\omega_{SRF}^D)^{n-m}}{\mid S\mid^2})$ and using $(*)$ we get the desired result.

\;

\;

\textbf{Remark B}:Note that the log semi-Ricci flat metric $\omega_{SRF}^D$ is not continuous in general. But if the central fiber has at worst log canonical singularities then, semi-Ricci flat metric is smooth in an open Zariski subset.

\textbf{Remark C}:So by applying the previous remark, the relative volume form $$\Omega_{(X,D)/Y}=\frac{(\omega_{SRF}^D)^n\wedge\pi^*\omega_{can}^{m}}{\pi^*\omega_{can}^m\mid S\mid^2}$$
is not smooth in general, where $S\in H^0(X,L_N)$ and $N$ is a divisor which come from canonical bundle formula of Fujino-Mori.

\;

\textbf{Remark D}: K\"ahler potential of Weil-Petersson metric induces a singular Hermitian metric with semi-positive curvature current on the tautological quotient bundle over the projective-space bundle $\mathbb P(f_*(K_{X/B}))$.

\;

Now we explain that under some algebraic condition the K\"ahler potential of Weil-Petersson metric on the moduli space of log Calabi-Yau pairs may be continuous (due to Tsuji-Kawamata). We recall the following Kawamata's theorem. \cite{19}
\begin{thm}
Let $f : X\to B$ be a surjective morphism of smooth projective
varieties with connected fibers. Let $P=\sum_jP_j$, $Q=\sum_lQ_l$, be normal crossing
divisors on $X$ and $B$, respectively, such that $f^{-1}(Q)\subset P$ and $f$ is smooth over $B\setminus Q$. Let $D =\sum_j d_jP_j$ be a $\mathbb Q$-divisor on $X$, where $d_j$ may be positive, zero or negative, which satisfies the following conditions A,B,C:

\;

\textbf{A}) $D = D^h + D^v$ such that any irreducible component of $D^h$ is mapped surjectively
onto $B$ by $f$ , $f : Supp(D^h) \to B$ is relatively normal crossing over $B \setminus Q$, and
$f(Supp(Dv ))\subset  Q$. An irreducible component of $D^h$ (resp. $D^v$ ) is called horizontal
(resp. vertical)

\;

\textbf{B})$d_j < 1$ for all $j$

\;

\textbf{C}) The natural homomorphism $\mathcal O_B \to f_*\mathcal O_X(\lceil -D\rceil)$
is surjective at the generic point of $B$.

\textbf{D}) $K_X + D\sim_{\mathbb Q} f^*
(K_B + L)$ for some $\mathbb Q$-divisor $L$ on $B$.

Let

\begin{align*}
f^*Q_l &= \sum_jw_{lj}P_j \\
 \bar{d_j} &=\frac{d_j+w_{lj}-1}{w_{lj}},\; \text{if}\;  f(P_j)=Q_l \\
\delta_l  &= max\{ \bar{d_j}; f(P_j)=Q_l \}.\\
\Delta  &= \sum_l \delta_l Q_l .\\
M  &= L-\Delta .\\
\end{align*}

Then $M$ is nef.
\end{thm}

The following theorem is straightforward from Kawamata's theorem

\begin{thm}Let $d_j < 1$ for all $j$ be as above in Theorem 11, and fibers be log Calabi-Yau pairs, then $$\int_{X_s\setminus D_s}(-1)^{n^2/2}\frac{\Omega_s\wedge\overline{\Omega_s}}{\mid S_s\mid^2}$$  is continuous on a nonempty Zariski open subset of $B$. 

\end{thm}

Since the inverse of volume gives a singular hermitian line bundle, we have the following theorem from Theorem 11

\begin{thm}Let $K_X + D\sim_{\mathbb Q} f^*
(K_B + L)$ for some $\mathbb Q$-divisor $L$ on $B$ and

\begin{align*}
f^*Q_l &= \sum_jw_{lj}P_j \\
 \bar{d_j} &=\frac{d_j+w_{lj}-1}{w_{lj}},\; \text{if}\;  f(P_j)=Q_l \\
\delta_l  &= max\{ \bar{d_j}; f(P_j)=Q_l \}.\\
\Delta  &= \sum_l \delta_l Q_l .\\
M  &= L-\Delta .\\
\end{align*}

Then 

$$\left(\int_{X_s\setminus D_s}(-1)^{n^2/2}\frac{\Omega_s\wedge\overline{\Omega_s}}{\mid S_s\mid^2}\right)^{-1}$$
is a continuous hermitian metric on the $\mathbb Q$-line bundle $K_B + \Delta$ when fibers are log Calabi-Yau pairs.

\end{thm}

\textbf{Remark E}: Note that Yoshikawa\cite{31}, showed that when the base of Calabi-Yau fibration $f:X\to B$ is a disc and central fibre $X_0$ is reduced and irreducible and pair $(X,X_0)$ has only canonical singularities then the K\"ahler potential of Weil-Petersson can be extended to a continuous Hermitian metric lying in the following class(due to Y.Kawamata also)

$$\mathcal B(B)=C^\infty(S)\oplus \bigoplus_{r\in \mathbb Q\cap (0,1]}\bigoplus_{k=0}^n\mid s\mid^{2r}(\log \mid s\mid)^kC^\infty(B)$$

In fact.
$$\int_{X_s\setminus D_s}(-1)^{n^2/2}\frac{\Omega_s\wedge\overline{\Omega_s}}{\mid S_s\mid^2}=C\left(\log\mid t\mid\right)^m\mid t\mid^{2k}(1+o(1))$$
where $C$ is a constant \cite{64}, \cite{57}, \cite{63}.

\;

If the central fiber be log Calabi-Yau pair, then $m=0$ and vise versa. Let $X\to \Delta$ be a one dimensional degeneration of Calabi-Yau varieties. If $X_0$ has finite Weil–Petersson distance, then after running a
log minimal model program, $(X_0, D_0)$ can be replaced by a log Calabi–Yau variety with at worst
log canonical singularities and vise versa after finite base change and birational change, see\cite{76}, \cite{69}. So this tells us that after running log MMP and semi-stable reduction, the logarithmic Weil-Petersson distance from $0$, has the following asymptotic formula, 

$$\int_{X_s\setminus D_s}(-1)^{n^2/2}\frac{\Omega_s\wedge\overline{\Omega_s}}{\mid S_s\mid^2}=C\mid t\mid^{2k}(1+o(1))$$

Note that, If $X_0$ only has canonical singularities, or if $X$
is smooth and $X_0$ only has isolated ordinary quadratic singularities, then
if $\pi:X\to \mathbb C^*$ be a family of degeneration of of Calabi-Yau fibers. Then the $L^2$-metric

$$\int_{X_s}\Omega_s\wedge\bar\Omega_s$$ is continuous. See Remark 2.10. of \cite{68}

Now, such fiberwise integral are not $C^\infty$ in general, and their smoothness correspond to Monodromy representation theory.
\;
\;

\textbf{Smoothness of fiberwise integral of Calabi-Yau volume}: Let $X$ be a closed normal analytic subspace in some open subset $U$ of $\mathbb C^N$ with an isolated singularity. Take $f:X\to \Delta$ be a degeneration of smooth Calabi-Yau manifolds, then $$s\to\int_{X_s}\Omega_s\wedge\bar\Omega_s\in C^\infty$$ if and only if the monodromy $M$ acting on the cohomology of the Milnor fibre of $f$ is the
identity and the restriction map  $j: H^n(X^*) \to H^n(F)^M$ is surjective, where $X^*=X\setminus \{0\}$ and $M$ denotes monodromy acting on $H^n(F)$ and $H^n(F)^M$ is the $M$-invariant subgroup and $F$ is the Milnor fiber at zero(see Corollary 6.2. \cite{72}). In fact the $C^\infty$-smoothness of fiberwise Calabi-Yau volume $\omega_{SRF}^m$ must correspond to such information of D.Barlet program.

\;

\textbf{Remark F}: Note that hermitian metric of semi Ricci flat metric $\omega_{SRF}$ is in the class of $\mathcal B(B)$ and is smooth if and only if we the central fiber has canonical singularities (due to diameter bound of fibers near to central fiber for obtaining $C^0$-estimate)

Now, we can give a parallel definition of Tsuji's foliation for fiberwise log Calabi-Yau metric.

\begin{defn} The null direction semi Ricci flat metric $\omega_{SRF}$ gives a foliation along Iitaka fibration $\pi:X\to Y$ and we call it fiberwise Calabi-Yau foliation and can be defined as follows $$\mathcal F=\{\theta\in T_{X/Y}|\omega_{SRF}(\theta,\bar\theta)=0\}$$
and along log Iitaka fibration  $\pi:(X,D)\to Y$, we can define the following foliation $$\mathcal F'=\{\theta\in T_{X'/Y}|\omega_{SRF}^D(\theta,\bar\theta)=0\}$$
where $X'=X\setminus D$. In fact from Theorem 0,9. the Weil-Petersson metric $\omega_{WP}$ vanishes everywhere if and only if $\mathcal F=T_{X/Y}$
\end{defn}

From Tsuji's work we have

\textbf{Lemma:} Let $\mathcal L$ be a leaf of $f_*\mathcal F'$, then $\mathcal L$ is a closed complex submanifold and the leaf $\mathcal L$ can be seen as fiber on the moduli map $$\eta:\mathcal Y\to \mathcal M_{CY}^D$$ where $\mathcal M_{CY}^D$ is the moduli space of log calabi-Yau fibers with at worst canonical singularites and 

$$\mathcal Y=\{y\in Y_{reg}|(X_y,D_y)\; \; \text{is Kawamata log terminal pair}\}$$

The following definition introduced by Tsuji 

\begin{defn}
Let $X$ be a compact complex manifold and let $L$ be a line bundle on $X$. A singular Hermitian metric $h$ on $L$ is said to be an analytic Zariski decomposition(or shortly AZD), if the following hold.

1. the curvature $\Theta_h$ is a closed positive current.
\;

2. for every $m\geq 0$, the natural inclusion 

$$H^0(X,\mathcal O_X(mL)\otimes \mathcal I(h^m))\to H^0(X,\mathcal O_X(mL))$$
is an isomorphism, where $\mathcal I(h^m)$ denotes the multiplier ideal sheaf of $h^m$.

\end{defn}

Since the Weil-Petersson metric is semi-positive.

\textbf{Remark G:} The hermitian metric corresponding to Song-Tian measure is Analytic Zariski Decomposition., i.e.,

$$h=\left(\frac{(\omega_{SRF}^D)^n\wedge\pi^*\omega_{can}^m}{\pi^*\omega_{can}^m\mid S\mid^2}\right)^{-1}$$

is AZD
\;

Now we show that Song-Tian-Tsuji measure is bounded if and only if central fiber of Iitaka fibration has canonical singularities at worst. Note that Song-Tian just showed that such measure is in $L^{1,\epsilon}$ and they couldn't prove the boundedness of such measure.

\begin{thm}Song-Tian-Tsuji measure 

$$\Omega_{X/Y}=\left(\frac{\omega_{SRF}^n \wedge\pi^*\omega_{can}^m}{\pi^*\omega_{can}^m}\right)$$

on Iitaka fobration is bounded if and only if the central fiber $X_0$ has log terminal singularities.
\;
\;

Proof: R. Berman in \cite{40}, showed that for Iitaka fibration $F:X\to Y$ a canonical sequence of Bergman type measures

$$v_k=\int_{X^{N_k-1}}\mu^{(N_k)}$$

where $$\mu^{(N_k)}=\frac{1}{Z_k}|S^{(k)}(z_1,...,z_{N_k})|^{2/k}dz_1\wedge d\overline {z}_1... dz_{N_k}\wedge d\overline z_{N_k}$$
and $N_{Z_k}$
is the normalizing constant ensuring that
$\mu^{(N_k)}$
is a
probability measure, and 

$$S^{(k)}(x_1,...,x_{N_k}):=\det \left(s_i^{(k)}(x_j)\right)_{1\leq i,j\leq N_k}$$ where $s_i^{(k)}$ is a basis in $H^0(X,kK_X)$

then $v_k$ converges weakly to Song-Tian measure $$\mu_X=F^*(\omega_Y)^{\kappa(X)}\wedge \omega_{SRF}^{n-\kappa(X)}$$

But from Proposition 1.17 \cite{27},( \cite{26} also), we know that If $K_{X_0}$ is a $\mathbb Q$–Cartier divisor then $X_0$ has log terminal singularities if and only if $\mathcal L^{2,m}_{X_0} = \mathcal O_{X_0}(m(K_{X_0}))$ for all
$m\geq 1$ where $L^{2,m}$ is the sheaf of locally $2/m$-integrable m-fold d-norms. So Song-Tian measure is bounded if and only if central fiber has log terminal singularities (in log Iitaka fibration case from Proposition 1.19 \cite{27},  we need to assume log canonical singularities for central fiber to obtain boundedness of Song-Tian measure).

\end{thm}

Mumford in his celebrated paper \cite{59} introduced the notion of good metric (which we call such metrics as Mumford metric) to extend the Chern-Weil theory to quasi-projective manifolds. First we recall Mumford metrics. For fixing our notation. Let $\bar E$
be a holomorphic vector bundle of rank $l$ over $X$ and $E=\bar E|_X$ and  $h$ an Hermitian metric on $E$ which may be singular near divisor $D$. Now cover a neighborhood of $D\subset X$ by finitely many polydiscs $\{U_\alpha=(\Delta^n,(z_1,...,z_n))\}$ such that $V_\alpha=U_\alpha\setminus D=(\Delta^*)^k\times \Delta^{n-k}$.  Namely, $U_\alpha\cap D=\{z_1...z_k=0\}$. Let $U=\cup_\alpha U_\alpha$ and $V=\cup V_\alpha$ then on each $V_\alpha$ we have the local Poincaré metric.

\begin{defn}
Let $\omega$ be a smooth local $p$-form defined on $V_\alpha$ we say
$\omega$ has Poincaré growth if there is a constant
$C_\alpha>0$ depending on $\omega$ such that

$$|\omega(t_1,t_2,...,t_p)|^2\leq C_\alpha \prod_{i=1}^p||t_i||^2$$
for any point
$z
\in V_\alpha$
and $t_1,...,t_p\in T_zX$ where $||.||$ is taken on Poincaré metric. We say $\omega$
is Mumford metric if both
$\omega$ and $d\omega$ have Poincaré growth.
\end{defn}
Now we recall the notion Hermitian Mumford  metric
$h$ on the vector bundle $E$.

\begin{defn}

An Hermitian metric
$h$
on
$E$
is good if for all
$z
\in
V$
, assuming
$z
\in
V_\alpha$
, and for all basis $(e_1,...,e_n)$ of $E$ over $U_\alpha$, if we let $h_{i\bar j}=h(e_i,e_{\bar j})$, then

\quad
\quad

I)$|h_{i\bar j}|, \left(\det h\right)^{-1}\leq C(\sum_{i=1}^k\log |z_i|)^{2n}$ for some $C>0$

\quad
\quad

II) The local 1-forms $(\partial h.h^{-1})_{\alpha\lambda}$ are Mumford on $V_\alpha$, namely the local 

connection
and curvature forms of $h$ have Poincaré growth.
\end{defn}

 If
$h$
is a good metric on
$E$
, the Chern forms $c_i(E,h)$ are
Mumford forms.

\;
\;

\begin{thm} Fiberwise Calabi-Yau metric on Iitaka fibration and hermitian metric corresponding to Song-Tian-Tsuji measure (inverse of Song-Tian-Tsuji measure $h=\Omega_{X/Y}^{-1}$ gives a singular hermitian metric) is a good metric in the sense of Mumford, i.e, it is Mumford metric when central fiber has log terminal singularities at worst.
\end{thm}

\;
Proof: From the Berman's formula in our proof of Theorem 16, and boundedness of Song-Tian-Tsuji measure and using Theorem 30, Lemma 36 in \cite{60}, we get the desired result.

\;
\;

\textbf{Remark}: Fiberwise log Calabi-Yau metric on log Iitaka fibration and hermitian metric corresponding to Song-Tian-Tsuji measure (inverse of Song-Tian-Tsuji measure $h=\Omega_{(X,D)/Y}^{-1}$ gives a singular hermitian metric) is a good metric in the sense of Mumford, i.e, it is Mumford metric when central fiber has log canonical singularities at worst.

\;

The following Theorem, gives a result to the goodness of canonical metric in the sense of Mumford, on moduli spaces of Calabi-Yau varieties with mild singularities in the sense of Minimal Model Prgram.

\begin{thm}
Since if a singular hermitian metric be a good metric in the sense of Mumford (Mumford metric), its first Chern current is also good metric in the sense of Mumford, hence from Theorem 9, the singular Weil-Petersson metric on moduli space of Calabi-Yau varieties is a good metric when central fiber has log terminal singularities at worst. Moreover if central fiber has log canonical singularities at worst, then the singular logarithmic Weil-Petersson metric on the moduli spaces of log Calabi-Yau varieties is a good metric in the sense of Mumford
\end{thm}
\textbf{$\bullet$ Weil-Petersson metric via Bismut-Vergne localization formula }

\;
\;

Let $\pi : M \to S$ be a submersion of smooth manifolds, with compact oriented fibres $X_s=\pi^{-1}(s)$. Let $T$ be a torus acting smoothly on $M$, and preserving the fibers $X_s$. Fix $K\in \mathfrak t$. Let $X_{s,K} \subset  X$ be the zero set of $K^{X_s}$ where $K^{X_s}$ is the corresponding
vector fields on $X_s$. Then $X_{s,K}$ is a totally geodesic submanifold
of $X_s$. 

 Let fibers $X_s$ have vanishing first Chern class $c_1(X_s)=0$. So we have the following formula of Weil-Petersson metric on moduli space of Calabi-Yau manifolds

$$\omega_{WP}=\sqrt{-1}\partial\bar\partial\log\int_{X_s}\Omega_s\wedge \overline {\Omega_s}$$

So, from Bismut-Vergne localization formula \cite{65},\cite{66}, \cite{67} we have the following formula for Weil-Petersson metric

$$\omega_{WP}=\sqrt{-1}\partial\bar\partial\log\int_{X_{s,K}}\frac{\Omega_s\wedge \overline {\Omega_s}}{e_K(N_{X_{s,K}/X_{s}}, \nabla^{N_{X_{s,K}/X_s}})}\;\; \text{mod}\; \; d^S\Omega(S)$$

where the equivariant Euler class is as follows

$$e_K(N_{X_{s,K}/X_{s}}, \nabla^{N_{X_{s,K}/X_s}})=\text{det}^{\frac{1}{2}}\left[\frac{J_K+R^{N_{X_{s,K}/X_s}}}{2\pi}\right]$$

where $J_K$ is an anti-symmetric parallel endomorphism
of $N_{X_{s,K}} /X_s$ , which is nondegenerate and $R^{N_{X_{s,K}/X_s}}$ is the curvature of any metric connection in $N_{X_{s,K}/X_s}$

\;
\;

In fact such K\"ahler potential of Weil-Petersson metric

$$\int_{X_{s,K}}\frac{\Omega_s\wedge \overline {\Omega_s}}{e_K(N_{X_{s,K}/X_{s}}, \nabla^{N_{X_{s,K}/X_s}})}\;\; \text{mod}\; \; d^S\Omega(S)$$

is important to rewrite the relative volume form $\Omega_{M/S}$  by using such K\"ahler potential to get canonical metric on total space along Iitaka fibration(when $\mathbb T$ be a torus acting smoothly on $M$, and preserving the fibers also) and reduce twisted version of Complex-Monge Ampere equation to modulo base space $S$ (along Iitaka fibration or canonical model $X\to X_{can}$).

\;
\;

Now, if we assume the fibers $(X_s, L_s)$ are polarized Calabi-Yau manifolds, then there exists $m_0 = m_0(n) > 0$, such
that for any $m \geq m_0$, then we can embed $i_s:X_s\to \mathbb P(H^0(X_s, -m(K_{X_s}+L_s))^*)$. Now we choose an orthonormal basis ${s_i(s)}$ of $ \mathbb P(H^0(X_s, -m(K_{X_s}+L_s))^*)$ with $L^2$-inner product, then we can write $$\Omega_{X_s}:=\Omega_s=\left(\sum_{i=1}^{N_m}|s_i(s)|^2\right)^{-\frac{1}{m}}$$

So we can rewrite the Weil-Petersson metric as

$$\omega_{WP}=\sqrt{-1}\partial\bar\partial\log\int_{X_{s,K}}\frac{\left(\sum_{i=1}^{N_m}|s_i(s)|^2\right)^{-\frac{1}{m}}\wedge \overline {\left(\sum_{i=1}^{N_m}|s_i(s)|^2\right)^{-\frac{1}{m}}}\otimes(ds\otimes d\bar s)^{-1}}{e_K(N_{X_{s,K}/X_{s}}, \nabla^{N_{X_{s,K}/X_s}})}\;\; \text{mod}\; \; d^S\Omega(S)$$

\section{Invariance of Plurigenera and positivity of logarithmic Weil-Petersson metric}

\;
Let $X$ be a projective variety, for every positive integer $m$, the $m$-th
plurigenus $P_m(X)$ is defined by

$$P_m(X)=\dim H^0(X,\mathcal O_X(mK_X))$$
The plurigenera are fundamental invariants under birational transformation of compact complex manifolds. It is an open conjecture that plurigenera is invariant under K\"ahler deformations.  For non-K\"ahler manifolds we don't have the invariance of plurigenera. For projective varieties this conjecture has been solved. 

We can solve Invariance of plurigenera by using algebraic geometric method and also by using complex analytical method. 

By using Tsuji-Boucksom's method \cite{50},\cite{51}, A complete solution for deformation of invariance of plurigenera is equivalent with the existence of
singular hermitian metric $h$ on $K_X$ such that the curvature current $\Theta_h$ is
semipositive and $h|_{X_t}$
is an AZD  of $K_{X_t}$
for every $t \in \Delta$.
As soon as we construct such a metric $h$, the $L^2$
-extension theorem implies the invariance of the plurigenera.

In algebraic language invariance of plurigenera is equivalent as follows.

Algebraic language of Invariance of plurigenera say that any section $mK_{X_0}$ extends to $X$, in other words, the restriction map 
$$H^0(X,mK_X)\to H^0(X_0,mK_{X_0})$$
is surjective.

Takayama, showed that, if $X_0$ has at most terminal singularities, then
$X_t$ has terminal singularities at most for all $t\in B$ . Moreover,If
$\pi:X\to B$
is smooth and the “abundance conjecture” holds true for general
$X_t$,then
$$P_m(X_t)=\dim H^0(X_t,mK_{X_t})$$ is independent of $t\in B$ for all $m$.
N.Nakayama, showed that the invariance of plurigenera for smooth projective deformations can be derived if the minimal model program were
completed for families immediately
follow from the $L^2$-extension theorem.

Now, we apply the relative K\"ahler-Ricci flow method which is the analytical version of MMP to solve invariance of plurigenera and give a method for the problem of deformation of invariance of plurigenera in K\"ahler setting.(see the preprint about Semipositivity of relative canonical bundles via K\"ahler-Ricci
flows, from Boucksom, S.; Tsuji, Hajime \cite{50})

Let
$\pi: X\to \Delta $ is smooth, and every
$X_t$ is of general type over open disc $\Delta$. Consider the relative K\"ahler Ricci flow 

$$\frac{\partial\omega}{\partial t}=-Ric_{X/\Delta}(\omega)-\Phi\omega$$

where $\Phi$ is the fiberwise constant function.

If such flow has solution, then we get the invariance of plurigenera. In fact if such flow has solution then $K_{X/\Delta}$ is psudo-effective. In fact if such flow has solution then $e^{-t}\omega_0+(1-e^{-t})c_1(K_{X/\Delta})$ is a K\"ahler metric and since we take the initial metric to be positive, hence $c_1(K_{X/\Delta})$ is positive. Now since fibers are of general type then by using Yau's theorem there exists a K\"ahler-Einstein metric on each fiber. So we can introduce fiberwise K\"ahler Einstein metric $\omega_{SKE}$ which by using Schumacher and Paun  result(a solution for such flow gives the semi-positivity of such fiberwise K\"ahler-Einstein metric. In fact to rewrite our relative K\"ahler Ricci flow to complex Monge-Ampere equation, we need to construct relative volume form by using such metric and relative volume form $\Omega_{X/\Delta}=\omega_{SKE}^m\wedge \pi^*(dt)$ must be semi-positive),  $\omega_{SKE}$ is semi-positive and its corresponding hermitian metric is AZD, hence we can take a hermitian metric $h^{\omega_{SKE}}$ on $K_X$ such that the curvature current $\Theta_{h^{\omega_{SKE}}}$ is
semipositive and $h^{\omega_{SKE}}|_{X_t}$
is an AZD  of $K_{X_t}$
for every $t \in \Delta$.
and by the $L^2$
-extension theorem(it is still open for K\"ahler manifolds) we get the invariance of the plurigenera and we obtain Theorem of Siu.

In singular setting, we need to the following conjecture

Let
$\pi: X\to \Delta$ is a family of projective varieties, and every
$X_t$ is of general type over open disc $\Delta$. Then fiberwise K\"ahler-Einstein metric is semi-positive if and only if the central fiber $X_0$ has at worst canonical singularities.

If this Conjecture mentioned before, holds true, then we can have invariance of plurigenera of Kawamata result in singular setting. In fact to obtain a $C^0$ solution for relative K\"ahler Ricci flow we must have diameter bound of fibers $X_t$ and this must be equivalent with canonical singularities of central fiber $X_0$. A weak answer of this is the theorem of Donaldson-Sun \cite{58}. A theorem of Donaldson-Sun states that if $X_t$ are K\"ahler-Einstein metric with negative Ricci curvature with uniform diameter bound, then the central fiber is normal and Kawamata log terminal singularities at worst. In view of the moduli theory of canonically polarized varieties, limit of fibers should have canonical singularities.

\;

\;

\;

\;

\end{document}